\documentclass[leqno,12pt]{amsart}
\usepackage{amscd,amsfonts,amsmath,amssymb,amstext,amsthm}
\newtheorem{theorem}{Theorem}[section]
\newtheorem{lemma}[theorem]{Lemma}
\newtheorem{proposition}[theorem]{Proposition}

\begin{document}
\title{An Estimate for the Average Number of Common Zeros of Laplacian Eigenfunctions }
\author{Dmitri Akhiezer and Boris Kazarnovskii}
\dedicatory{To Ernest Borisovich Vinberg on the occasion of his 80th birthday}
\address {Institute for Information Transmission Problems \newline
19 B.Karetny per.,127994, Moscow, Russia,\newline
{\rm D.A.:} 
{\it akhiezer@iitp.ru},  
{\rm B.K.:} 
{\it kazbori@gmail.com}.}
\thanks{The research was carried out at the Institute for Information Transmission Problems under support by 
the Russian Foundation of Sciences, 
grant  No. 14-50-00150}
\begin{abstract}
On a compact Riemannian manifold $M$ of dimension $n$, we consider $n$ eigenfunctions of the Laplace operator $\Delta $ with eigenvalue $\lambda$. If $M$ is homogeneous under a compact Lie group preserving the metric then 
we prove that the average number of common zeros of $n$ eigenfunctions does not exceed 
$c(n)\lambda^{n/2}{\rm vol}\,M$, the expression known from the celebrated Weyl's law.
Moreover, if the isotropy representation is irreducible, then the estimate turns into the equality.
The constant  $c(n)$ is
explicitly given.
The method of  proof is based on the application of Crofton's formula for the sphere. \end{abstract}
\renewcommand{\subjclassname}
{\textup{2010} Mathematics Subject Classification}
\maketitle
\renewcommand{\thefootnote}{}
%\footnotetext{Supported by }
\maketitle
\section{Introduction}\label{intro}
Let $M$ be a compact Riemannian manifold 
without boundary, $n = {\rm dim}\, M$, and $dx$ the Riemannian measure on $M$. For an eigenvalue
$\lambda $ of the Laplace operator $\Delta $ on $M$  let $H(\lambda )$ 
denote the corresponding eigenspace, i.e.,
$$H(\lambda) =\{ f\in C^\infty (M, \mathbb R)\ \vert \ \Delta u + \lambda u = 0\}.$$
Then $H(\lambda)$ is a finite dimensional real vector subspace of $L^2(M, dx)$, 
considered with the induced scalar product. We note that the space $H(\lambda)$ and the scalar product  are 
invariant under any isometry of $M$.
Our goal is to define and, under certain assumptions, to evaluate the average number
of zeros of the system of equations
$$u_1 =u_2 =\ldots=u_n =0,$$
where $u_i \in H(\lambda )$
are linearly independent.
The linear envelope of $u_i$ is a subspace
$U\subset H(\lambda )$  of dimension $n$. We define $Z(U)$ 
to be the set of common
zeros of $u_i$.
 For any set $A$ we denote by $\#A$ the number of its points if $A$ is finite or $\infty $ if $A$ is infinite. 
We now let $U$ vary within a subspace $H \subset H(\lambda)$ of some dimension $N>0$,
which does not necessarily coincide with $H(\lambda)$. Then we have the function
$$U \mapsto \#Z(U) \in \mathbb N \cup \{\infty \}.$$ The average number of zeros ${\mathfrak M}\{\#Z(U)\} $ or, more precisely,   
${\mathfrak M}_H\{\#Z(U)\}   $
is defined as the integral of $\#Z(U)$ over the Grassmanian ${\rm Gr}_n(H)$ with respect to the measure
induced by the (normalized) Haar measure of ${\rm SO}(N, \mathbb R)$ acting on $H$. We keep in mind that the averaging process
depends on $H$, but do not stress this in notation. 

\begin{theorem}\label{main1} Let $M=K/V$ be a homogeneous space of a compact connected Lie group $K$
with a $K$-invariant Riemannian metric and let $H \subset H(\lambda )$ be a $K$-invariant subspace.
 Then
$$\mathfrak{M}\{\#Z(U)\} \le c(n) {\lambda}^{n/2}
{\rm vol}\,M,$$
where 
$$c(n) ={ 2\over{\sigma_nn^{n/2}}}$$
and $\sigma _n$ is the volume of the $n$-dimensional sphere of radius 1.
Furthermore, if  points of $M$ are not locally separated by $H$ then
$$\mathfrak {M}\{ \#Z(U)\} = 0.$$

\end{theorem}

Recall that $M=K/V$ is called isotropy irreducible if the representation of the isotropy subgroup $V$
in the tangent space at the origin $e\cdot V$ is irreducible. We refer the reader to \cite{GHL} for the
properties of isotropy irreducible homogeneous spaces. Their minimal immersions into spheres are for the first time
considered in \cite{Ta}. We also note that a symmetric space of a simple group, e.g.,
the sphere with the special orthogonal group, is isotropy irreducible. Under
a stronger assumption, namely, if the connected component  $V_0\subset V$ acts irreducibly in the tangent
space at $e\cdot V$,
all isotropy irreducible
homogeneous spaces are listed, see \cite{{Ma},{Wo}}.

\begin{theorem}\label{main2} If $M=K/V$ is isotropy irreducible then for any $K$-invariant subspace $H \subset H(\lambda)$ one has the equality
$$\mathfrak{M}\{\#Z(U)\} = c(n) {\lambda}^{n/2}
{\rm vol}\,M.$$
\end{theorem}
This is proved in \cite{AK}, but our proof here is slightly different. Namely, Theorem \ref {main2}
appears below as a corollary of Theorem \ref{main1}. 
It should be noted that the inequality in Theorem \ref{main1} is in general strict even if $H$ locally separates points of $M$,
see an example in Sect.\ref{conc}.

As in \cite{AK}, our arguments are based
on Crofton's formula for the sphere \cite{Sa}, which we recall in
Sect.\ref{Cr}.
In integral geometry, formulae of Crofton type are regarded as the simplest kinematic formulae,
see \cite{{Ho},{Sa}}. We refer the reader to \cite{{GS},{PF}} for another approach,
namely, for the relationship of Crofton's
formula with Radon transform.

Finally, we want to mention two circumstances which partly motivated our research.
We recall that the zero set of an eigenfunction of the Laplace operator is called
a nodal set. Connected components of its complement are called nodal domains.
The classical Courant's theorem \cite{CH} says that the number of nodal domains defined by
the $k$-th eigenfunction does not exceed $k$, see also \cite {Ch}. Consider now the
set $Z$ of common zeros of $m\le n$ eigenfunctions. As a rule, its complement for $m>1$ is connected.
In order to carry over Courant's theorem to this case,
V.Arnold suggested to study the topology of the analytic set $Z$ 
and to find the dependence of suitable topological invariants of $Z$ 
on the number of the corresponding eigenvalue of the Laplacian, see \cite{Ar}, Problem 2003-10, p.174. 
We follow Arnold's suggestion for $m=n$ under certain additional assumptions of group-theoretic character.
On the other hand,
we remark that the right hand side of the inequality in Theorem \ref{main1}, up to a coefficient, coincides
with the first term of asymptotics in the celebrated Weyl's law, see \cite{Iv}.
Thus we get an estimate of the average number of common zeros of Laplacian eigenfunctions in terms 
of the asymptotic expression for the eigenvalue number. There is no doubt that one has to give a meaning to this fact.

\section{ Crofton's formula for the sphere }\label{Cr}

Let $S = S^{N-1} \subset {\mathbb R}^N$ be the unit sphere considered with the metric induced from the ambient
Euclidean space.
Suppose we are given two submanifolds $M\subset S$ and $L\subset S$ of dimensions $n$ and $l$
respectively, such that $n+l = N-1$. Let $\sigma _n$ denote the volume of the $n$-dimensional sphere 
of radius 1. Moving the second submanifold by a rotation $r \in {\rm SO}(N, \mathbb R)$,
take the number of intersection points ${\#}(M \cap r
\cdot L)$. Then the Crofton (or kinematic) formula is the expression 
for the average number of such points, i.e., the integral
of ${\#}(M\cap r\cdot L)$ over ${\rm SO}(N,\mathbb R)$ with respect to the Haar measure $dr$.
Namely,
$$\int_{M\cap r\cdot L \ne \emptyset}\ {\#}(M \cap r\cdot L)\cdot dr = {2\over \sigma _n \sigma_ l}\cdot ({\rm vol}\, M)({\rm vol}\, L).$$
We refer the reader to the book of L.A.Santal\'o \cite{Sa}, see Sections 15.2 and 18.6. Another proof is given by R.Howard, see
\cite{Ho}, Sect. 3.12. We will use the Crofton formula for $L$ being a plane section through the origin of complementary
dimension $l=N-n-1$. Then the formula reads
$$\int_{M\cap r\cdot L\ne \emptyset }\ {\#}(M \cap r\cdot L) \cdot dr = {2\over \sigma_n}\cdot {\rm vol}\,M.$$
It will be convenient to have this formula for coverings.
For this we recall the definition of $k$-density on a manifold $X$,
see, e.g.,\,\cite{PF}.  Let $\Lambda_X^k(x) $ be the cone of decomposable $k$-vectors
at $x \in X$. Then $k$-density on $X$ is a function $\mu: \Lambda_X^k(x) \to {\mathbb R}$,  depending smoothly on $x$
and such that $\mu (t\xi) = \vert t \vert \mu (\xi)$. A $k$-density can be integrated along any, not necessarily orientable,
submanifold of dimension $k$. If $X$ is a Riemannian manifold then for any $k$ there is a $k$-density 
attaching to a $k$-vector $\xi_1\wedge \cdots \wedge \xi_k$ the
volume of a $k$-dimensional parallelotope with edges $\xi_1, \ldots, \xi_k$. For any differentiable mapping $\phi: Y \to X$ and for a
$k$-density $\mu $ on $X$ we have the inverse image $\phi^*
\mu $, a $k$-density on $Y$.

\begin{proposition}\label{cr} Let $M$ be a compact manifold of dimension $n$ and let $\phi : M \to S = S^{N-1} \subset {\mathbb R}^N$
be a covering over $\phi (M)$. 
Assume that $\phi (M)$ is not contained in a proper linear subspace of ${\mathbb R}^N$.
The
pull-back of the dual space $({\mathbb R}^N){^*}$
is some $N$-dimensional subspace
$H \subset C^\infty (M, \mathbb R)$, which comes equipped with the metric of $({\mathbb R}^N){^*}$.
 Then, in the notation of Sect.\ref{intro},
for an $n$-dimensional subspace $U\subset H$ one has
$${\mathfrak M}\{\#Z(U)\}= {2\over\sigma_n}\int_M\phi ^*\Omega_n,$$
where $\Omega _n$ is the Riemannian $n$-density on $S$.
\end{proposition}

\medskip\noindent
{\it Proof.} We have
$$Z(U) = \phi ^{-1}(\phi(M)\cap L)$$
for some plane section $L$ of $S$ of codimension $n$. Let $d$ be the degree of the covering $\phi: M \to \phi (M)$. Then
$${\mathfrak M}\{\#Z(U)\} = {2d\over \sigma_n} \cdot {\rm vol}\,\phi(M) = {2d\over \sigma_n} \int_{\phi(M)}\Omega _n
={2\over\sigma_n}\int _M\phi^*\Omega _n$$
by Crofton's formula and by the definition of $\Omega _n$.
\hfill{$\square $}

\section{Equivariant mappings into the sphere}\label{equi}

Let $K$ be a compact connected Lie group , $V \subset K$ a closed subgroup, and
$M =K/V$ the homogeneous space considered with the transitive $K$-action. Introduce a $K$-invariant
Riemannian metric $g$ on $M$. Let $H$ be a $K$-invariant finite dimensional
subspace in $ C^\infty (M, \mathbb R)$ considered with the scalar product
$$(f_1,f_2) = \int_Mf_1(x)f_2(x) \cdot dx,$$
where $dx$ is the Riemannian measure on $M$. For $x\in M$ let $\alpha _x$ denote the linear functional from $H^*$
attaching to a function $f\in H$ its value $f(x)$. Take an orthonormal basis $\{f_i\}_{i=1}^N$ of $H$,
denote by $\{f_i^*\}$ the dual basis of $H^*$,
 and identify
$H^*$ with ${\mathbb R}^N$ via
$$(\mu_1, \ldots, \mu_N)\mapsto \mu_1f_1^* + \ldots + \mu_Nf_N^*.$$
Then 
$$\alpha _x = \sum _{i=1}^N f_i(x)f_i^*,$$
and so we obtain the coordinate form of the map $M\ni x \to \alpha_x\in H^*$, namely,
$$\Phi=(f_1,\ldots,f _ N): M \to {\mathbb R}^N.$$
The map $\Phi$ is equivariant with respect to the given action of $K$ on $M$ and the linear
representation of $K$ in $H^*$ identified with ${\mathbb R}^N$. In fact, the image of $\Phi $
is contained in a sphere, and we now compute its radius. The following lemma 
for spherical harmonics is known in quantum mechanics as Uns\"old's identity (1927).
The general case is found in \cite{GHL}, Exercise 5.25 c),i), p.261 and p.303, see also \cite {AK}.

\begin{lemma}\label {rad}The image $\Phi(M)$ is contained in the sphere
of radius $R = \sqrt{N\over {\rm vol}\, M}$, i.e.,
$$ \sum_{i=1}^N\, f_i^2 = {N\over{\rm vol}\,M}\,.$$ 
\end{lemma}
\noindent {\it Proof.} Write
$$f_i(k^{-1}x) = \sum _j\, a_{ij}(k)f_j(x).$$  
Then $(a_{ij})$ is an orthogonal matrix,
hence
$$\sum_i f_i^2(k^{-1}x) =
= \sum _if_i^2(x).$$
Since $K$ acts transitively on $M$, we have
$$\sum _if_i^2 = C,$$ 
where $C$ is some constant. Integrating the last equality over $M$, we get $N = C{\rm vol}\,M$.
\hfill {$\square $} 

\medskip
\noindent
The following lemma is quite general and elementary. For two quadratic forms $P$ and $Q$ on a real vector space and
$Q$ positive definite,
let ${\rm tr}_Q\,P$ denote the trace of the matrix of $P$ in any orthonormal basis with respect to $Q$.
We use the same notation for a quadratic and the associate symmetric bilinear form.
We denote by $g_0$ the Eucledian metric on $\mathbb R^N$, i.e., 
$g_0 =\sum dt_i^2$ relative to Cartesian coordinates.
\begin{lemma}\label{trace}Let $(M,g)$ be a Riemannian manifold, $$\Phi = (f_1,\ldots,f_N):M \to {\mathbb R}^N$$
a differential mapping, and $\Phi^*g_0 $ the pull-back of the Euclidean metric onto $M$. Then
$${\rm tr}_g\, \Phi^*(g_0) = \sum _{k=1}^N \ \vert \vert {\rm grad}\,f_k\vert \vert ^2,$$
where the norm of a tangent vector is taken with respect to $g$.
\end{lemma} 
\noindent{\it Proof.}\ Let $\xi _1, \ldots, \xi_n$ be an orthonormal basis at some point $x \in M$. Then 
$$h_{kj} = df_k(\xi_j) = g({\rm grad} f_k, \xi_j),\ \ j=1,\ldots,n,$$ 
are the coordinates of ${\rm grad} f_k$ relative to the chosen basis.
At the point $x$ we have
$$\Phi ^*(g_0)\,(\xi_i ,\xi _j) = g_0(d\Phi(\xi _i), d\Phi (\xi_j))= \sum_{k=1}^Ndf_k(\xi_i)df_k(\xi_j)=
\sum_{k=1}^Nh_{ki}h_{kj},$$ hence 
$${\rm tr}_g\, \Phi^*(g_0) =
\sum_{i=1}^n\sum_{k=1}^N\, h_{ki}^2=\sum _{k=1}^N\sum_{i=1}^n h_{ki}^2=
 \sum _{k=1}^N \ \vert \vert {\rm grad}\,f_k\vert \vert ^2$$
as we asserted.
\hfill $\square$

\medskip
\noindent
We now return to the Lie group setting.
Let $M = K/V$, where $K$ is a compact connected Lie group and $V \subset K$
a closed subgroup. Fix a $K$-invariant Riemannian metric $g$ on $M$, consider the
corresponding Laplace operator and take a $K$-invariant subspace $H \subset H(\lambda)$.
Recall that we have the mapping $\Phi: M \to {\mathbb R}^N$ and the quadratic form
 $\Phi ^*(g_0)$ on $M$ induced by the Euclidean metric $g_0$ on ${\mathbb R}^N$.
The first and the last assertions in the following lemma are known, see \cite{Ta}.

\begin{lemma}\label{trhom} {\rm (1)} The form $\Phi ^*(g_0)$ is $K$-invariant.\\ 
{\rm (2)} One has
$${\rm tr}_g\Phi ^*(g_0) = {\lambda N \over {\rm vol}\, M}\,.$$
{\rm (3)} If $M$ is isotropy irreducible then $\Phi^*(g_0) = c\cdot g$, where $c$ is a constant. 
\end{lemma}
\noindent{\it Proof.} 
(1)  follows from the fact that $\Phi $ is a $K$-equivariant mapping.
To prove (3), it suffices to remark that an irreducible $V$-module
has only one, up to a multiple, $V$-invariant quadratic form. 
To prove (2), write
$${\rm div}(f_i{\rm grad} \, f_i)= f_i\Delta f_i + \vert\vert{\rm grad}\, f_i\vert \vert^2 = -\lambda f_i^2 + 
\vert \vert {\rm grad}\, f_i \vert \vert ^2,$$
sum up the equalities for all $i$
and integrate over $M$. This gives
$$\int_M\sum_{i=1}^N\, \vert \vert {\rm grad}\, f_i\vert \vert ^2 dx = \lambda \sum_{i=1}^N \int_M f_i^2\,dx = \lambda N.$$
On the other hand,
${\rm tr}_g\Phi^*(g_0)$ is constant, and so we obtain
$$\int_M\sum_{i=1}^N\, \vert \vert {\rm grad}\, f_i\vert \vert ^2 dx = \int_M {\rm tr}_g\Phi ^*(g_0)\cdot dx = 
{\rm tr}_g\Phi ^*(g_0)\cdot {\rm vol}\, M$$
by Lemma \ref{trace}. Therefore
${\rm tr}_g\Phi^*(g_0) 
\cdot {\rm vol}\, M = \lambda N$.
\hfill{$\square$}

\medskip
\noindent
{\it Remark.}\ If $M$ is isotropy irreducible then $\Phi $ is a minimal
immersion of $M$ into the sphere, see \cite{Ta}.

\section{Proof of main results} 
\medskip
\noindent{\it Proof of Theorem \ref{main1}.}

\noindent
Case 1: {\it $H$ locally separates the points of $M$.} The equivariant mapping $\Phi $ defined in Sect.\ref{equi} gives rise to a
fibering $\Phi : M\to \Phi(M)$. By Lemma \ref{rad} the image is contained in the
sphere of radius $R = \sqrt{N\over {\rm vol}\, M}$. Therefore $x \mapsto \phi (x) = \sqrt{{ {\rm vol}\,M \over N}}\Phi (x)$
is a covering over a submanifold $\phi (M)\subset S $. By Proposition \ref{cr}
$$\mathfrak{M}\{\#Z(U)\}  = {2\over \sigma _n} \int_M\phi ^*\Omega _n.$$
Now, $\phi ^*\Omega _n$ can be obtained by pulling the Riemannian metric $g_0$ back to $M$ by $\phi $ and
then taking the volume form. Suppose $\phi^ *(g_0)$ has eigenvalues $\beta _1, \ldots, \beta_n$
with respect to an orthogonal basis of $g$ on $M$ at some point $x \in M$. 
Then $\phi ^*\Omega _n = \sqrt{\beta_1\cdots\beta_n}\cdot dx$.
By (1) of Lemma \ref{trhom} the form $\phi ^*(g_0)$
is $K$-invariant, and so the set of eigenvalues (with multiplicities) is the same for all $x\in M$.
Thus
$${\mathfrak M}\{\#Z(U)\} = {2\over \sigma_n} \sqrt{\beta _1 \cdots \beta_n}\cdot {\rm vol}\,M.$$ 
From (2) of Lemma \ref{trhom} it follows that
$\beta _1 + \ldots +\beta _n = \lambda$.
Since
$$(\beta _1 \cdots  \beta _n)^{1/n} \le {{\beta _1 + \ldots + \beta _n}\over n},$$
we obtain the estimate
$$\mathfrak M \{\#Z(U)\} \le {2\over \sigma _n}\Bigl ({\lambda \over n}\Bigr)^{n/2}\cdot{\rm vol}\, M$$
as required.

\noindent
Case 2: {\it The points of $M$ are not locally separated by $H$.} The equivariant mapping $\Phi $ has positive
dimensional fiber, i.e., $$m := {\rm dim}\, \Phi(M) < {\rm dim}\, M = n.$$
Recall that $\Phi (M) \subset S \subset H^*$, ${\rm dim}\, H = N$, and $l=N-1-n$.
Consider the mapping of
$$\Gamma =\{(x,L)\,\vert \,x \in \Phi (M), L \in {\rm Gr}_{l+1}(H^*), x \in L\}$$
to ${\rm Gr}_{l+1}(H^*) \simeq {\rm Gr}_n(H)$
defined by $(x,L) \mapsto L$. 
It is easily seen that $\Gamma $ is a manifold. Namely, $\Gamma $ is the space of a fiber bundle
over $\Phi (M)$ whose fiber is the set of points of the Grassmanian ${\rm Gr}_{l+1}(H^*)$,
such that the corresponding subspace contains a given point $x \in \Phi(M)$. Thus the
dimension of $\Gamma $  is easily computed, namely,
$${\rm dim} \, \Gamma = {\rm dim}\, \Phi(M) + {\rm dim}\, {\rm Gr}_l(H_0^*),$$
where $H_0^*$
is the quotient space of $H^*$ of dimension $N-1$. Hence
$${\rm dim}\, \Gamma = m +(N-1-l)l < n + n(N-1-n) = {\dim \rm Gr}_n(H).$$
Therefore the image of $\Gamma $ in ${\rm Gr}_n(H)$ has measure 0 (for example, by Sards's lemma), and
our assertion follows by the definition of ${\mathfrak M}\{\#Z(U)\}$.

\hfill$\square$

\bigskip
\noindent
{\it Proof of Theorem \ref{main2}}. Since $M= K/V$ is isotropy irreducible,  $V$ is not contained in a proper subgroup of $K$ of
greater dimension. Thus the mapping $\Phi $ is a covering, so that we are in Case 1 of Theorem \ref{main1}.
Also, (3) of Lemma \ref{trhom} shows that all $\beta _i$ are equal.
Thus the inequality in the proof of Theorem \ref{main1}, Case 1,   becomes an equality.
\hfill{$\square$}
\section{Concluding remarks}\label{conc}
\noindent 1.\,{\it Example: the inequality in Theorem \ref{main1} can be
strict.}
 Let $e_1, e_2$ be an orthonormal basis in ${\mathbb R}^2$, $x,y$ the Cartesian coordinates relative to $e_1, e_2$,
and $\Gamma \subset {\mathbb R}^2$ the lattice generated by $2\pi e_1$ and $(2\pi/a)e_2$, where $a $ is some fixed non-zero number.   
The metric $dx^2 + dy^2$ defines a flat metric on the torus $M = {\mathbb R}^2/\Gamma$.
Let $\Delta $ be the associate Laplace operator. Then the minimal
positive eigenvalue of $\Delta $ equals $\lambda _1 =1+a^2$,
the space $H(\lambda _1)$ has dimension 4, and one finds easily an orthonormal basis $\{f_i\}$ in $H(\lambda _1)$. Namely,
$$f_1 = {\sqrt{a}\over \pi}\, {\rm sin}\,x
\cdot {\rm cos}\, ay,\  f_2 = {\sqrt{a}\over \pi}\, {\rm sin}\,x
\cdot {\rm sin}\, ay,$$
$$  f_3 = {\sqrt {a}\over \pi}\, {\rm cos}\,x
\cdot {\rm cos} \,ay, \ f_4 = {\sqrt {a}\over \pi}\, {\rm cos}\,x
\cdot {\rm sin}\, ay.$$
 Then
$$df_1^2 + df_2^2+df_3^2+df_4^2 ={a\over \pi^2} (dx^2 + a^2dy^2)$$
by a direct computation.
The mapping $\phi: M \to S^3 \subset {\mathbb R}^4$ from the proof of Theorem \ref{main1}
has the form 
$$(x,y) \mapsto {\pi \over \sqrt{a}}( f_1(x,y), f_2 (x,y), f_3 (x,y), f_4(x,y)).$$ 
Using the notation inroduced in the course of the proof, we get
$\beta _1 =1, \beta _2 =a^2$, so the inequality is strict if $a\ne 1$.

\medskip
\noindent
2. {\it Crofton's formula is local.} In Proposition \ref{cr}, one can replace $M$ by any domain $D\subset M$ and get the
average number
${\mathfrak M}^D\{\#Z(U)\}$ of common zeros of $u_i$ in $D$. Namely, 
$$\mathfrak{M}^D\{\#Z(U)\} =
{2\over\sigma_n}\int_D\phi ^*
\Omega _n.$$  For $M=K/V$ this reduces to
$$\mathfrak{M}^D\{\#Z(U)\} =
{2\over \sigma _n}\sqrt{\beta_1\cdots \beta_n}\cdot {\rm vol}\, D ,$$ see the proof
of Theorem \ref{main1}. Finally, if $M$ is isotropy irreducible,
then
$$\mathfrak{M}^D\{\#Z(U)\} = {2\over\sigma_n}\Bigl({\lambda\over n}\Bigr)^{n/2}\cdot {\rm vol}\, D.$$

\medskip
\noindent
3. {\it The case of $k$ equations, $ k < n$.} We keep the notations and assumptions  of Theorem \ref{main2}, but
consider a smaller number $k$ of functions $u_i$. 
Let  the functions $u_i \in H \subset H(\lambda),\ i = 1,2,\ldots,k,$
be linearly independent so that their linear envelope $U\subset H(\lambda)$ has
dimension $k$. We keep the same notation $Z(U)$ for the set of common zeros of $u_i$.
If $Z(U)$ is a smooth $(n-k)$-dimensional submanifold of $M$ then we denote by ${\rm vol}_{n-k}\, Z(U)$ its
$(n-k)$-dimensional volume. Otherwise we put ${\rm vol}_{n-k}\, Z(U) = 0$. 
The subset of ${\rm Gr}_k(H)$ 
consisting of all $U$, such that $Z(U)$ is not a smooth submanifold of dimension $n-k$ has measure 0. 
This follows from Sard's lemma applied to the projection
$$\Gamma = \{(x,U)\, \vert \, x\in Z(U)\} \to {\rm Gr}_k(H),\ \ (x,U) \mapsto U.$$ 
Therefore the average volume of $Z(U)$ is correctly defined as the integral of  ${\rm vol}_{n-k}\, Z(U)$ over
the Grassmanian ${\rm Gr}_k(H)$. For $M=K/V$ isotropy irreducible we have 
$${\mathfrak M}\{{\rm vol}_{n-k}\, Z(U)\} = {\sigma_{n-k}\over \sigma_n}\Bigl( {\lambda\over n}\Bigr)^{k/2} \cdot {\rm vol}\, M.$$   
The proof follows the above scheme for $n$ functions if one takes into account the following fact.
Let $M$ be an $n$-dimensional submanifold of the sphere $S\subset {
\mathbb R}^N$ and let $L\subset {\mathbb R}^N$ 
be a vector subspace of codimension $k$. Then it follows from Crofton's formula
that the integral of ${\rm vol}_{n-k} \, M\cap L$ over the Grassmanian
of all such subspaces equals
${\sigma _{n-k}\over \sigma _n}\cdot {\rm vol}\, M$.

\bigskip\noindent{\it Note.} V.M.Gichev informed the authors
that Theorem \ref{main2} follows from the results of \cite{Gi}.

\end {document}